\documentclass[journal]{IEEEtran}
\usepackage[dvips]{graphicx}
\usepackage[font=footnotesize,labelfont=bf]{caption}
\usepackage{amsmath,amssymb,bm}
\usepackage{algorithm,algorithmic}
\usepackage{amsfonts,dsfont,color,bbm}
\usepackage{epstopdf,epsf,epsfig}
\usepackage{cite,subfig}
\usepackage[percent]{overpic}

\renewcommand{\vec}[1]{\mbox{\boldmath${#1}$}}
\newcommand{\norm}[1]{\left|\left|#1\right|\right|}
\newcommand{\nn}{\nonumber}
\newcommand{\Real}{{\mathbbm{R}}}
\newcommand{\E}{{\mathbbm{E}}}

\newcommand{\W}{\mathcal{W}}

\DeclareMathOperator*{\argmin}{arg\,min}

\newcommand{\vx}{\vec{x}}
\newcommand{\vg}{\vec{g}}
\newcommand{\vbg}{\vec{\bar{g}}}
\newcommand{\vw}{\vec{w}}
\newcommand{\vu}{\vec{u}}
\newcommand{\vbw}{\vec{\bar{w}}}
\newcommand{\vpsi}{\vec{\psi}}
\newcommand{\vphi}{\vec{\phi}}
\newcommand{\vH}{\vec{H}}
\newcommand{\vX}{\vec{X}}
\newcommand{\vG}{\vec{G}}
\newcommand{\vW}{\vec{W}}

\newcommand{\vA}{\vec{A}}
\newcommand{\vB}{\vec{B}}
\newcommand{\vC}{\vec{H}}

\begin{document}

\title{Stochastic Subgradient Algorithms for Strongly Convex Optimization over Distributed Networks}

\author{N. Denizcan Vanli, Muhammed O. Sayin, and Suleyman S. Kozat
\thanks{This work is supported in part by Outstanding Researcher Program of Turkish Academy of Sciences and in part by TUBITAK under Contract 113E517.}
\thanks{N. D. Vanli is with the Laboratory for Information and Decision Systems, Massachusetts Institute of Technology, Cambridge, MA 02139 (e-mail: denizcan@mit.edu).}
\thanks{M. O. Sayin is with the Coordinated Science Laboratory, University of Illinois at Urbana-Champaign, Urbana, IL 61801 (e-mail: sayin2@illinois.edu).}
\thanks{S. S. Kozat is with the Department of Electrical and Electronics Engineering, Bilkent University, Ankara 06800, Turkey (e-mail: kozat@ee.bilkent.edu.tr).} }

\maketitle

\begin{abstract}
We study diffusion and consensus based optimization of a sum of unknown convex objective functions over distributed networks. The only access to these functions is through stochastic gradient oracles, each of which is only available at a different node, and a limited number of gradient oracle calls is allowed at each node. In this framework, we introduce a convex optimization algorithm based on the stochastic gradient descent (SGD) updates. Particularly, we use a carefully designed time-dependent weighted averaging of the SGD iterates, which yields a convergence rate of $O\left(\frac{N\sqrt{N}}{T}\right)$ after $T$ gradient updates for each node on a network of $N$ nodes. We then show that after $T$ gradient oracle calls, the average SGD iterate achieves a mean square deviation (MSD) of $O\left(\frac{\sqrt{N}}{T}\right)$. This rate of convergence is optimal as it matches the performance lower bound up to constant terms. Similar to the SGD algorithm, the computational complexity of the proposed algorithm also scales linearly with the dimensionality of the data. Furthermore, the communication load of the proposed method is the same as the communication load of the SGD algorithm. Thus, the proposed algorithm is highly efficient in terms of complexity and communication load. We illustrate the merits of the algorithm with respect to the state-of-art methods over benchmark real life data sets and widely studied network topologies.
\end{abstract}
\begin{keywords}
Distributed processing, convex optimization, online learning, diffusion strategies, consensus strategies.
\end{keywords}

\section{Introduction}\label{sec:introduction}
\IEEEPARstart{T}{he} demand for large-scale networks consisting of multiple agents (i.e., nodes) with different objectives is steadily growing due to their increased efficiency and scalability compared to centralized distributed structures \cite{sayed_mag1,sayed_mag2,giannakis_mag,regret}. A wide range of problems in the context of distributed and parallel processing can be considered as a minimization of a sum of objective functions, where each function (or information on each function) is available only to a single agent or node \cite{rel1,rel2,rel3}. In such practical applications, it is essential to process the information in a decentralized manner since transferring the objective functions as well as the entire resources (e.g., data) may not be feasible or possible \cite{sayed1,sayed2,sayed3,sayed4}. For example, in a distributed data mining scenario, privacy considerations may prohibit sharing of the objective functions \cite{rel1,rel2,rel3}. Similarly, in a distributed wireless network, energy considerations may limit the communication rate between agents \cite{sayed5,sayed6,sayed7,sayed8}. In such settings, parallel or distributed processing algorithms, where each node performs its own processing and share information subsequently, are preferable over the centralized methods \cite{sayed9,sayed10,sayed11,sayed12}.

Here, we consider minimization of a sum of unknown convex objective functions, where each agent (or node) observes only its particular objective function via the stochastic gradient oracles. Particularly, we seek to minimize this sum of functions with a limited number of gradient oracle calls at each agent. In this framework, we introduce a distributed online convex optimization algorithm based on the SGD iterates that efficiently minimizes this cost function. Specifically, each agent uses a time-dependent weighted combination of the SGD iterates and achieves the presented performance guarantees, which matches with the lower bounds presented in \cite{agarwal}, only with a relatively small excess term caused by the unknown network model. The proposed method is comprehensive, in that any communication strategy, such as the diffusion \cite{sayed_mag1} and the consensus \cite{regret} strategies, are incorporated into our algorithm in a straightforward manner as shown in the paper. We compare the performance of our algorithm respect to the state-of-the-art methods \cite{regret,sayed2,rakhlin} in the literature and present the outstanding performance improvements for various well-known network topologies and benchmark data sets.

Distributed networks are successfully used in wireless sensor networks \cite{giannakis1,giannakis2,giannakis3,giannakis4,giannakis5,giannakis6}, and recently used for convex optimization via projected subgradient techniques \cite{regret,sayed1,sayed2,rel1,rel2,rel3}. In \cite{sayed2}, the authors illustrate the performance of the least mean squares (LMS) algorithm over distributed networks using different diffusion strategies. We emphasize that this problem can also be casted as a distributed convex optimization problem, hence our results can be applied to these problems in a straightforward manner. In \cite{sayed1}, the authors consider the cooperative optimization of the cost function under convex inequality constraints. However, the problem formulation as well as the convergence results in this paper are significantly different than the ones in \cite{sayed1}. In \cite{regret}, the authors present a deterministic analysis of the SGD iterates and our results builds on them by illustrating a stronger convergence bound in expectation while also providing MSD analyses of the SGD iterates. In \cite{rel1,rel2,rel3}, the authors consider the distributed convex optimization problem and present the probability-$1$ and mean square convergence results of the SGD iterates. In this paper, on the other hand, we provide the expected convergence rate of our algorithm and the MSD of the SGD iterates at any time instant.

Similar convergence analyses are recently illustrated in the computational learning theory \cite{hazan1,hazan2,agarwal,rakhlin,lacoste}. In \cite{hazan1}, the authors provide deterministic bounds on the learning performance (i.e., regret) of the SGD algorithm. In \cite{hazan2}, these analyses are extended and a regret-optimal learning algorithm is proposed. In the same lines, in \cite{rakhlin}, the authors describe a method to make the SGD algorithm optimal for the strongly convex optimization. However, these approaches rely on the smoothness of the optimization problem. In \cite{lacoste}, a different method to achieve the optimal convergence rate is proposed and its performance is analyzed. On the other hand, in this paper, the convex optimization is performed over a network of localized learners, unlike \cite{hazan1,hazan2,rakhlin,lacoste}. Our results illustrate the convergence rates over any unknown communication graph, and in this sense build upon the analyses of the centralized learners. Furthermore, unlike \cite{hazan2,rakhlin}, our algorithm does not require the optimization problem to be sufficiently smooth.

Distibuted convex optimization appears in a wide range of practical applications in wireless sensor networks and real-time control systems \cite{sayed_mag1,sayed_mag2,giannakis_mag}. We introduce a comprehensive approach to this setup by proposing an online algorithm, whose expected performance is asymptotically the same as the performance of the optimal centralized processor. Our results are generic for any probability distribution on the data, not necessarily Gaussian unlike the conventional works in the literature \cite{sayed2,sayed3}. Furthermore, our performance bounds are optimal in a strong deterministic sense up to constant terms. Our experiments over different network topologies, various data sets and cost functions illustrate the superiority and robustness of our approach with respect to the state-of-the-art methods in the literature.

Our main contributions are as follows.
\begin{enumerate}
  \item We introduce a distributed online convex optimization algorithm based on the SGD iterates that achieves an optimal convergence rate of $O\left(\frac{N\sqrt{N}}{T}\right)$ after $T$ gradient updates, for each and every node on the network. We emphasize that this convergence rate is optimal since it achieves the lower bounds presented in \cite{agarwal} up to constant terms.
  \item We show that the average SGD iterate achieves a MSD of $O\left(\frac{\sqrt{N}}{T}\right)$ after $T$ gradient updates.
  \item Our analyses can be extended to analyze the performances of the diffusion and consensus strategies in a straightforward manner as illustrated in our paper.
  \item We illustrate the highly significant performance gains of the introduced algorithm with respect to the state-of-the-art methods in the literature under various network topologies and benchmark data sets.
\end{enumerate}

The organization of the paper is as follows. In Section \ref{sec:prob}, we introduce the distributed convex optimization framework and provide the notations. We then introduce the main result of the paper, i.e., a SGD based convex optimization algorithm, in Section \ref{sec:main} and analyze the convergence rate of the introduced algorithm. In Section \ref{sec:sim}, we demonstrate the performance of our algorithm with respect to the state-of-the-art methods through simulations and then conclude the paper with several remarks in Section \ref{sec:conc}.

\section{Problem Description}\label{sec:prob}

\subsection{Notation}
Throughout the paper, all vectors are column vectors and represented by boldface lowercase letters. Matrices are represented by boldface uppercase letters. For a matrix $\vH$, $\norm{\vH}_F$ is the Frobenius norm. For a vector $\vx$, $\norm{\vx} = \sqrt{\vx^T \vx}$ is the $\ell^2$-norm. Here, $\vec{0}$ (and $\vec{1}$) denotes a vector with all zero (and one) elements and the dimensions can be understood from the context. For a matrix $\vH$, $\vH_{ij}$ represents its entry at the $i$th row and $j$th column. For a convex set $\W$, $\Pi_\W$ denotes the Euclidean projection onto $\W$, i.e., $\Pi_\W(\vw') = \argmin_{\vw\in\W} \norm{\vw-\vw'}$.

\subsection{System Overview}
We consider the problem of distributed strongly convex optimization over a network of $N$ nodes. At each time $t$, each node $i$ observes a pair of regressor vectors and data, i.e., $\vu_{t,i}$ and $d_{t,i}$, where the pairs $(\vu_{t,i}, d_{t,i})$ are independent and identically distributed (i.i.d.) for all $t\geq1$, which is a common framework in the conventional studies in the literature \cite{sayed2,sayed3,giannakis4}. Here, the distributions of the data pairs can differ from one node to another and we do not have any information on the underlying distributions. In such a framework, the aim of each node is to minimize a strongly convex cost function $f_i(\cdot)$ over a convex set $\W$. However, $f_i(\cdot)$'s are now known and each node accesses to its $f_i(\cdot)$ only via a stochastic gradient oracle, which given some $\vw\in\W$, produces a random vector $\hat{\vg}_i$, whose expectation $\E\hat{\vg}_i=\vg_i$ is a subgradient of $f_i$ at $\vw$. Using these stochastic gradient oracles, each node estimates a parameter of interest $\vw_{t,i}$ on a common convex set $\W$ and calculates an estimate of the output as
\begin{equation}
  \hat{d}_{t,i} = \vw_{t,i}^T \vu_{t,i}, \nn
\end{equation}
i.e., by a first order linear method. After observing the true data at time $t$, node $i$ suffers a loss of $f_{t,i}(\vw_{t,i})$, where $f_{t,i}(\cdot)$ is a strongly convex loss function. As an example, our cost function can be defined as follows
\begin{equation}\label{eq:ex}
  f_{t,i}(\vw_{t,i}) = \ell(\vw_{t,i};\vu_{t,i},d_{t,i}) + \frac{\lambda}{2} \norm{\vw_{t,i}}^2,
\end{equation}
where $\ell(\vw_{t,i};\vu_i,d_i)$ is a Lipschitz-continuous convex loss function with respect to the first variable, where \eqref{eq:ex} is extensively studied in the literature \cite{hazan1,hazan2,lacoste,rakhlin} as a strongly convex loss function involving regularity terms.

\begin{figure}
  \centering
  \includegraphics[width=.48\textwidth]{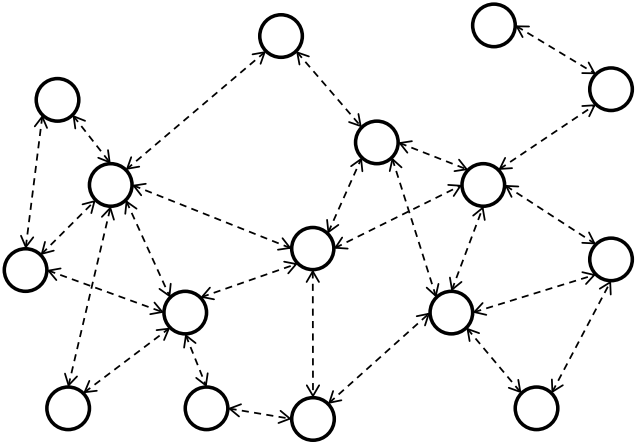}\\
  \caption{An example distributed network of $16$ nodes, where each agent communicates with a set of nodes in its neighborhood.}\label{fig:network}
\end{figure}

Here, the aim of each node is to minimize its expected loss over the convex set $\W$. To continue with our example in \eqref{eq:ex}, the aim of each node is to minimize
\begin{equation}\label{eq:loss}
  \E f_i(\vw_{t,i}) = \E \, \ell(\vw_{t,i};\vu_i,d_i) + \frac{\lambda}{2} \norm{\vw_{t,i}}^2,
\end{equation}
where we drop the notational dependency of $f$ to the time index $t$ since the data pairs are assumed to be independent over time. We emphasize that the formulation in \eqref{eq:loss} covers a wide range of practical loss functions. As an example, when $\ell(\vw_{t,i};\vu_i,d_i) = (d_i-\vw_{t,i}^T\vu_i)^2$, we consider the regularized squared error loss and when $\ell(\vw_{t,i};\vu_i,d_i) = \max\{0,1-d_i\vw_{t,i}^T\vu_i\}$, we consider the hinge loss. Since we make no assumptions on the loss function $f_i(\vw_{t,i})$ other than strong convexity, one can also use different loss functions with their corresponding gradients and our results still hold.

In this framework, at each time $t$, after observing the true data $d_{t,i}$, each node exchanges (i.e., diffuses) information with its neighbors in order to produce the next estimate of the parameter vector $\vw_{t+1,i}$. For example, in Fig. \ref{fig:network}, we have a distributed network consisting of $N=16$ nodes, where each node communicates with its neighbors. Although the aim of each node is to minimize its expected loss in \eqref{eq:loss}, the ultimate goal of the distributed network is to minimize the total expected loss, i.e.,
\begin{equation}\label{eq:cost}
  \E f(\vw_{t,i}) \triangleq \sum_{j=1}^N \E f_j(\vw_{t,i}),
\end{equation}
for each node $i$, with at most $t$ stochastic gradient oracle calls at each agent, where $N$ is the number of nodes in the distributed network. Here, we emphasize that the expected total loss in \eqref{eq:cost} is defined with respect to the data statistics at every node, whereas the parameter vector is trained using the gradient oracle calls at $i$th node and the information diffusion over the distributed network. Thus, our aim is to minimize a global, strongly convex cost function at each node over a distributed network, where each learner is allowed to use at most $t$ calls to the gradient oracle until time $t$.

\section{Main Results}\label{sec:main}
In this section, we present the main result of the paper, where we introduce an algorithm based on the SGD updates that achieves an expected regret upper bound of $O\left(\frac{N\sqrt{N}}{T}\right)$ after $T$ iterates. The proposed method uses time dependent weighted averages of the SGD updates at each node together with the adapt-then-combine diffusion strategy \cite{sayed2} to achieve this performance. However, as later explained in the paper (i.e., in Section \ref{sec:ext}), our algorithm can be extended the consensus strategy in a straightforward manner. The complete description of the algorithm can be found in Algorithm \ref{alg}.

\begin{algorithm}[t]
    \caption{Time Variable Weighting (TVW)}\label{alg}
    \begin{algorithmic}[1]
        \FOR{$t=1$ \TO $T$}
            \FOR{$i=1$ \TO $N$}
                \STATE $\hat{d}_{t,i} = \vbw_{t,i}^T \vu_{t,i}$ \% Estimation
                \STATE $\vg_{t,i} = \nabla_{\vw_{t,i}} f_i(\vw_{t,i})$ \% Gradient oracle call
                \STATE $\vpsi_{t+1,i} = \vw_{t,i} - \mu_t \vg_{t,i}$ \% SGD update
                \STATE $\vphi_{t+1,i} = \Pi_\W(\vpsi_{t+1,i})$ \% Projection
                \STATE $\vw_{t+1,i} = \sum_{j=1}^N \vH_{ji} \vphi_{t+1,j}$ \% Diffusion
                \STATE $\vbw_{t+1,i} = \frac{t}{t+2} \vbw_{t,i} + \frac{2}{t+2} \vw_{t,i}$ \% Weighting
            \ENDFOR
        \ENDFOR
    \end{algorithmic}
\end{algorithm}

To achieve the aforementioned result, we first introduce the following lemma, which provides an upper bound on the performance of the average parameter vector.\\

{\bf Lemma 1}
{\em Assume that (i) each $f_i$ is $\lambda$-strongly convex over $\W$, where $\W \subseteq \Real^p$ is a closed convex set; (ii) the norms of any two subgradients has a finite covariance, i.e., $\E\left[ \norm{\vg_i(\vw_1)}\norm{\vg_j(\vw_2)} \right] \leq G^2$, for any $\vw_1,\vw_2 \in \W$. Then defining
\begin{equation}
  \vw_t \triangleq \frac{1}{N} \sum_{i=1}^N \vw_{t,i}, \nn
\end{equation}
Algorithm \ref{alg} yields
\begin{align}\label{eq:lemma1}
  & \E\norm{\vw_{t+1}-\vw^*}^2 - (1-\lambda\mu_t)\,\E\norm{\vw_t-\vw^*}^2 \nn\\
    & \hspace{.5cm} \leq \frac{2\mu_t}{N} \, \E \Big[ f(\vw^*) - f(\vw_t) + G \sum_{i=1}^N \big( 2\norm{\vw_t-\vw_{t,i}}  \nn\\
    & \hspace{1.5cm} + \norm{\vw_t-\vpsi_{t+1,i}} \big) \Big] + 4G^2\mu_t^2.
\end{align}\\
}

This lemma provides an upper bound on the regret and the squared deviation of the average parameter vector. It provides an intermediate step to relate the performance of the parameter vector at each agent to the best parameter vector. The assumptions in Lemma 1 are widely used to analyze the convergence of online algorithms such as in \cite{hazan1,hazan2,rakhlin,lacoste}. Particularly, the first assumption in Lemma 1 is satisfied for many loss functions, such as the hinge loss, that are widely used in the literature. Even though this assumption may not hold for some loss functions such as the square loss, one can incorporate a small regularization term in order to make it strongly convex \cite{rakhlin,hazan1}. The second assumption in Lemma 2 is practically a boundedness condition that is widely used to analyze the performance of SGD based algorithms \cite{lacoste,hazan2}. We emphasize that our algorithm does not need to know this upper bound and it is only used in our theoretical derivations.\\

{\bf Proof }
In order to efficiently manage the recursions, we first consider the projection operation and let
\begin{equation}
  \vx_{t,i} \triangleq \Pi_\W(\vpsi_{t+1,i}) - \vpsi_{t+1,i},
\end{equation}
and $\vg_{t,i} \triangleq \nabla_{\vw_{t,i}} f_i(\vw_{t,i})$. Then, we can compactly represent the averaged estimation parameter $\vw_t = \frac{1}{N} \sum_{i=1}^N \vw_{t,i}$ in a recursive manner as follows \cite{regret}
\begin{align}
  \vw_{t+1} & = \frac{1}{N} \sum_{j=1}^N \left[ \sum_{i=1}^N \vH_{ij} \left( \vw_{t,i} - \mu_t \vg_{t,i} + \vx_{t,i} \right) \right] \nn\\ 
    & = \vw_t + \frac{1}{N} \sum_{i=1}^N \left( \vx_{t,i} - \mu_t \vg_{t,i} \right),
\end{align}
where the last line follows since $\vH$ is right stochastic, i.e., $\vH\vec{1}=\vec{1}$.

Hence, the squared deviation of these average iterates with respect to $\vw^*$ can be obtained as follows
\begin{align}\label{eq:1}
  \norm{\vw_{t+1}-\vw^*}^2 & =  \norm{\vw_t-\vw^* + \frac{1}{N} \sum_{i=1}^N \left( \vx_{t,i} - \mu_t \vg_{t,i} \right)}^2 \nn\\
    & \hspace{-1cm} = \norm{\vw_t-\vw^*} + \frac{1}{N^2} \norm{\sum_{i=1}^N (\vx_{t,i} - \mu_t \vg_{t,i})}^2 \nn\\
    & +\frac{2}{N} \sum_{i=1}^N (\vx_{t,i} - \mu_t \vg_{t,i})^T (\vw_t-\vw^*).
\end{align}

We first upper bound the second term in the right hand side (RHS) of \eqref{eq:1} as follows
\begin{align}\label{eq:bound1_old}
  & \frac{1}{N^2} \norm{\sum_{i=1}^N (\vx_{t,i} - \mu_t \vg_{t,i})}^2  \nn\\
    & \hspace{2cm} \leq \frac{1}{N^2} \left(\sum_{i=1}^N  \norm{\vx_{t,i}} + \mu_t \norm{\vg_{t,i}} \right)^2.
\end{align}
We then note that
\begin{align}
  \norm{\vx_{t,i}} & = \norm{\Pi_\W(\vpsi_{t+1,i}) - \vpsi_{t+1,i}} \nn\\
    & \leq \norm{\vw_{t,i} - \vpsi_{t+1,i}} \nn\\
    & = \mu_t \norm{\vg_t}, \nn
\end{align}
where the second line follows since
\begin{equation}
  \Pi_\W(\vpsi_{t+1,i}) = \argmin_{\vphi\in\W} \norm{\vphi-\vpsi_{t+1,i}}. \nn
\end{equation}
Thus, we can rewrite \eqref{eq:bound1_old} as follows
\begin{align}\label{eq:bound1}
  \frac{1}{N^2} \norm{\sum_{i=1}^N (\vx_{t,i} - \mu_t \vg_{t,i})}^2 & \leq \frac{4\mu_t^2}{N^2} \left(\sum_{i=1}^N  \norm{\vg_{t,i}} \right)^2  \nn\\
    & \leq 4 G^2 \mu_t^2,
\end{align}
where the last line follows since $\E\left[ \norm{\vg_{t,i}}\norm{\vg_{t,j}} \right] \leq G^2$ for any $i,j\in\{1,\dots,N\}$ according to the assumption.

We next turn our attention to $\left[-\vg_{t,i}^T(\vw_t-\vw^*)\right]$ term in \eqref{eq:1} and upper bound this term as follows
\begin{align}
  & -\vg_{t,i}^T(\vw_t-\vw^*) = - \vg_{t,i}^T(\vw_t-\vw_{t,i}+\vw_{t,i}-\vw^*) \nn\\
    & \hspace{.4cm} \leq - \vg_{t,i}^T(\vw_t-\vw_{t,i}) + f_i(\vw^*) - f_i(\vw_{t,i}) \nn\\
    & \hspace{1cm} - \frac{\lambda}{2} \norm{\vw^* - \vw_{t,i}}^2 \label{eq:2}\\
    & \hspace{.4cm} \leq - \vg_{t,i}^T(\vw_t-\vw_{t,i}) + \vbg_{t,i}^T(\vw_t-\vw_{t,i}) + f_i(\vw^*) \nn\\
    & \hspace{1cm} - f_i(\vw_t) - \frac{\lambda}{2} \norm{\vw^* - \vw_{t,i}}^2 - \frac{\lambda}{2} \norm{\vw_{t,i} - \vw_t}^2 \label{eq:3}\\
    & \hspace{.4cm} \leq f_i(\vw^*) - f_i(\vw_t) + \left( \norm{\vbg_{t,i}} + \norm{\vg_{t,i}} \right) \norm{\vw_t-\vw_{t,i}} \nn\\
    & \hspace{1cm} - \frac{\lambda}{2} \norm{\vw^* - \vw_{t,i}}^2 - \frac{\lambda}{2} \norm{\vw_{t,i} - \vw_t}^2, \label{eq:4}
\end{align}
where $\vbg_{t,i} \triangleq \nabla_{\vw_t} f_i(\vw_t)$, \eqref{eq:2} follows from the $\lambda$-strong convexity, i.e.,
\begin{equation}
  f_i(\vw^*) \geq f_i(\vw_{t,i}) + \vg_{t,i}^T (\vw^*-\vw_{t,i}) + \frac{\lambda}{2} \norm{\vw^* - \vw_{t,i}}^2, \nn
\end{equation}
\eqref{eq:3} also follows from the $\lambda$-strong convexity, i.e.,
\begin{equation}
  f_i(\vw_{t,i}) \geq f_i(\vw_t) + \vbg_{t,i}^T(\vw_{t,i}-\vw_t) + \frac{\lambda}{2} \norm{\vw_{t,i} - \vw_t}^2, \nn
\end{equation}
and \eqref{eq:4} follows from the Cauchy-Schwarz inequality. Summing \eqref{eq:4} from $i=1$ to $N$ and taking expectation of both sides, we obtain
\begin{align}\label{eq:bound2}
  & - \E \sum_{i=1}^N \vg_{t,i}^T(\vw_t-\vw^*) \nn\\
    & \hspace{.6cm} \leq \E \left[ 2G \sum_{i=1}^N \norm{\vw_t-\vw_{t,i}} + f(\vw^*) - f(\vw_t) \right. \nn\\
    & \hspace{1.2cm} \left. - \frac{\lambda N}{2} \sum_{i=1}^N \frac{1}{N} \left( \norm{\vw^* - \vw_{t,i}}^2 + \norm{\vw_{t,i} - \vw_t}^2 \right) \right] \nn\\
    & \hspace{.6cm} \leq \E \Big[ f(\vw^*) - f(\vw_t) + 2G \sum_{i=1}^N \norm{\vw_t-\vw_{t,i}} \nn\\
    & \hspace{1.2cm} - \frac{\lambda N}{2} \norm{\vw_t-\vw^*}^2 \Big],
\end{align}
where the first inequality follows due to the assumption and the last inequality follows from the Jensen's inequality due to the convexity of the norm operator.

We finally turn our attention to $\vx_{t,i}^T(\vw_t-\vw^*)$ term in \eqref{eq:1} and write this term as follows
\begin{align}
  \vx_{t,i}^T(\vw_t-\vw^*) & = \vx_{t,i}^T(\vw_t-\vpsi_{t+1,i}) + \vx_{t,i}^T(\vpsi_{t+1,i}-\vw^*) \nn\\
    & \leq \vx_{t,i}^T(\vw_t-\vpsi_{t+1,i}), \nn
\end{align}
where the inequality follows from the definition of the Euclidean projection. Taking the expectation of both sides, we can upper bound this term as follows
\begin{align}\label{eq:bound3}
  \E\,\vx_{t,i}^T(\vw_t-\vw^*) & = \E\left[\norm{\vx_{t,i}}\norm{\vw_t-\vpsi_{t+1,i}}\right] \nn\\
    & \leq G \mu_t \, \E \norm{\vw_t-\vpsi_{t+1,i}}.
\end{align}

Putting \eqref{eq:bound1}, \eqref{eq:bound2}, and \eqref{eq:bound3} back in \eqref{eq:1}, we obtain
\begin{align}\label{eq:step1}
  & \E\norm{\vw_{t+1}-\vw^*}^2 - (1-\lambda\mu_t)\,\E\norm{\vw_t-\vw^*}^2 \nn\\
    & \hspace{.5cm} \leq \frac{2\mu_t}{N} \, \E \Big[ f(\vw^*) - f(\vw_t) + G \sum_{i=1}^N \big( 2\norm{\vw_t-\vw_{t,i}}  \nn\\
    & \hspace{1.5cm} + \norm{\vw_t-\vpsi_{t+1,i}} \big) \Big] + 4G^2\mu_t^2.
\end{align}
This concludes the proof of Lemma 1. \hfill $\square$\\

Having obtained an upper bound on the performance of the average parameter vector, we then consider the MSD of the parameter vectors at each node from the average parameter vector. This lemma will then be used to relate the performance of each individual node to the performance of the fully connected distributed system.\\

{\bf Lemma 2}
{\em In addition to the assumption in Lemma 1, assume that (i) the communication graph $\vH$ forms a doubly stochastic matrix such that $\vH$ is irreducible and aperiodic; (ii) the initial weights at each node are identically initialized to avoid any bias, i.e., $\vw_{1,i}=\vw_{1,j}$, $\forall i,j \in \{1,\dots,N\}$. Then Algorithm \ref{alg} yields
\begin{equation}\label{eq:lemma2_1}
  \E\norm{\vw_t-\vw_{t,i}} \leq 2 G \sqrt{N} \sum_{z=1}^{t-1} \mu_{t-z} \sigma^z,
\end{equation}
and
\begin{equation}\label{eq:lemma2_2}
  \E\norm{\vw_t-\vpsi_{t+1,i}} \leq G\mu_t + 2 G \sqrt{N} \sum_{z=1}^{t-1} \mu_{t-z} \sigma^z,
\end{equation}
where $\sigma$ is the second largest singular value of matrix $\vH$.\\
}

We emphasize that the assumptions in Lemma 2 are not restrictive and previous analyses in the literature also use similar assumptions \cite{regret,rel1,rel2,rel3}. The first assumption in Lemma 2 indicates that the convergence of the communication graph is geometric. This assumption holds when the communication graph is strongly connected, doubly stochastic, and each node gives a nonzero weight to the iterates of its neighbors. This holds for many communication strategies such as the Metropolis rule \cite{sayed_mag1}. The second assumption in Lemma 2 is basically an unbiasedness condition, which is reasonable since the objective weight $\vw^*$ is completely unknown to us. Even though the initial weights are not identical, our analyses still hold, however with small additional excess terms.\\

{\bf Proof }
We first let $\vW_t \triangleq [\vw_{t,1},\dots,\vw_{t,N}]$, $\vG_t \triangleq [\vg_{t,1},\dots,\vg_{t,N}]$, and $\vX_t \triangleq [\vx_{t,1},\dots,\vx_{t,N}]$. Then, we obtain the recursion on $\vW_t$ as follows
\begin{equation}\label{eq:recursion}
  \vW_t = \vW_1 \vH^{t-1} + \sum_{z=1}^{t-1} \left( \vX_{t-z} - \mu_{t-z} \vG_{t-z} \right) \vH^z.
\end{equation}
Letting $\vec{e}_i$ denote the basis function for the $i$th dimension, i.e., only the $i$th entry of $\vec{e}_i$ is $1$ whereas the rest is $0$, we have
\begin{align}\label{eq:bound4}
  & \norm{\vw_t-\vw_{t,i}} = \norm{\vW_t\left( \frac{1}{N}\vec{1}-\vec{e}_i \right)} \nn\\
    & \hspace{.3cm} \leq \norm{\vw_1-\vw_{1,i}} + 2 \sum_{z=1}^{t-1} \mu_{t-z} \norm{\vG_{t-z}\left( \frac{1}{N}\vec{1}-\vH^{z}\vec{e}_i \right)} \nn\\
    & \hspace{.3cm} = 2 \sum_{z=1}^{t-1} \mu_{t-z} \norm{\vG_{t-z}\left( \frac{1}{N}\vec{1}-\vH^{z}\vec{e}_i \right)} \nn\\
    & \hspace{.3cm} \leq 2 \sum_{z=1}^{t-1} \mu_{t-z} \norm{\vG_{t-z}}_F \norm{\frac{1}{N}\vec{1}-\vH^{z}\vec{e}_i},
\end{align}
where the third line follows due to the unbiased initialization assumption, i.e., $\vw_{1,i}=\vw_{1,j}$, $\forall i,j\in\{1,\dots,N\}$.

We first consider the term $\norm{\frac{1}{N}\vec{1}-\vH^{z}\vec{e}_i}$ of \eqref{eq:bound4} and define the matrix $\vB \triangleq \frac{1}{K} \vec{1} \vec{1}^T$. Then, we can write
\begin{align}\label{eq:pol5}
  \norm{\frac{1}{N} \vec{1} - \vC^z \vec{e}_i} & = \norm{\vB\vec{e}_i-\vC^z\vec{e}_i} \nn\\
    & = \norm{(\vB-\vC)^z\vec{e}_i},
\end{align}
where the last line follows since $\vB^z=\vB$, $\forall z\geq1$.

Here, we let $\sigma_1(\vA) \geq \sigma_2(\vA) \geq \cdots \geq \sigma_N(\vA)$ denote the singular values of a matrix $\vA$ and $\lambda_1(\vA) \geq \lambda_2(\vA) \geq \cdots \geq \lambda_N(\vA)$ denote the eigenvalues of a symmetric matrix $\vA$. Then, we can upper bound \eqref{eq:pol5} as follows
\begin{equation}
  \norm{(\vB-\vC)^z\vec{e}_i} \leq \sigma_1(\vB-\vC) \norm{(\vB-\vC)^{z-1}\vec{e}_i}, \nn
\end{equation}
$\forall z\geq1$. Therefore, using the above recursion $z$ times to \eqref{eq:pol5}, we obtain
\begin{align}\label{eq:pol7}
  \norm{(\vB-\vC)^z\vec{e}_i} & \leq \sigma^z_1(\vB-\vC) \norm{\vec{e}_i} \nn\\
    & = \sigma^z_1(\vB-\vC).
\end{align}
Here, we note that $\vC$ is assumed to be a doubly stochastic matrix, thus $\lambda_1(\vC)=1$. Therefore, the eigenspectrums of $\vB-\vC$ and $(\vB-\vC)^T($\vB-\vC$)$ are equal to the eigenspectrums of $\vC$ and $\vC^T\vC$, respectively, except the largest eigenvalues, i.e., $\lambda_1(\vC) = \lambda_1(\vC^T\vC) = 1$. Thus, we have
\begin{equation}\label{eq:pol8}
  \sigma^z_1(\vB-\vC) = \sigma^z_2(\vC),
\end{equation}
and combining \eqref{eq:pol5}, \eqref{eq:pol7}, and \eqref{eq:pol8}, we obtain
\begin{equation}\label{eq:pol9}
  \norm{\frac{1}{N} \vec{1} - \vC^z \vec{e}_i} \leq \sigma^z_2(\vC).
\end{equation}
From here on, we denote $\sigma \triangleq \sigma_2(\vC)$ for notational simplicity.

Putting \eqref{eq:pol9} back in \eqref{eq:bound4}, we obtain
\begin{equation}\label{eq:bound6}
  \norm{\vw_t-\vw_{t,i}} \leq 2 \sum_{z=1}^{t-1} \mu_{t-z} \sigma^z \norm{\vG_{t-z}}_F.
\end{equation}
Taking the expectation of both sides and noting that
\begin{align}
  \left( \E\norm{\vG_{t-z}}_F \right)^2 & \leq \E\norm{\vG_{t-z}}_F^2 \nn\\
    & = \E\sum_{i=1}^N \norm{\vg_{t-z,i}}^2 \nn\\
    & \leq G^2N, \nn
\end{align}
we can rewrite \eqref{eq:bound6} as follows
\begin{equation}
  \E\norm{\vw_t-\vw_{t,i}} \leq 2 G \sqrt{N} \sum_{z=1}^{t-1} \mu_{t-z} \sigma^z.
\end{equation}

An upper bound for the term $\norm{\vw_t-\vpsi_{t+1,i}}$ can be obtained as
\begin{align}
  \norm{\vw_t-\vpsi_{t+1,i}} & = \norm{\vw_t-\vw_{t,i}+\mu_t\vg_{t,i}} \nn\\
    & \leq \norm{\vw_t-\vw_{t,i}} + \mu_t\norm{\vg_{t,i}}, \nn
\end{align}
where the last line follows from the triangle inequality. Taking expectation of both sides, we obtain the following upper bound
\begin{equation}\label{eq:bound5}
  \E\norm{\vw_t-\vpsi_{t+1,i}} \leq G\mu_t + 2 G \sqrt{N} \sum_{z=1}^{t-1} \mu_{t-z} \sigma^z.
\end{equation}
This concludes the proof of Lemma 2. \hfill $\square$\\

The results in Lemma 1 and Lemma 2 are combined in the following theorem to obtain a regret bound on the performance of the proposed algorithm. This theorem illustrates the convergence rate of our algorithm (i.e., Algorithm \ref{alg}) over distributed networks. The upper bound on the regret $O\left(\frac{N\sqrt{N}}{T}\right)$ results since the algorithm suffers from a ``diffusion regret'' to sufficiently exchange (i.e., diffuse) the information among the nodes. This convergence rate matches with the lower bounds presented in \cite{agarwal} up to constant terms, hence is optimal in a minimax sense. The computational complexity of the introduced algorithm is on the order of the computational complexity of the SGD iterates up to constant terms. Furthermore, the communication load of the proposed method is the same as the communication load of the SGD algorithm. On the other hand, by using a time-dependent averaging of the SGD iterates, our algorithm achieves a significantly improved performance as shown in Theorem 1 and illustrated through our simulations in Section \ref{sec:sim}.\\

{\bf Theorem 1}
{\em Under the assumptions in Lemma 1 and Lemma 2, Algorithm \ref{alg} with learning rate $\mu_t=\frac{2}{\lambda(t+1)}$ and weighted parameters $\vbw_{t,i}$, when applied to any independent and identically distributed regressor vectors and data, i.e., $(\vu_{t,i},d_{t,i})$ for all $t\geq1$ and $i=1,\dots,N$, achieves the following convergence guarantee
\begin{equation}\label{eq:theorem1}
  \E\left[ f\left( \vbw_{T,i} \right) - f(\vw^*) \right] \leq \frac{4NG^2}{\lambda(T+1)} \left( 3 + \frac{8\sigma\sqrt{N}}{1-\sigma} \right),
\end{equation}
for all $T \geq 1$, where $\sigma$ is the second largest singular value of matrix $\vH$.\\
}

This theorem illustrates that although every node uses local gradient oracle calls to train its parameter vector, each node asymptotically achieves the performance of the centralized processor through the information diffusion over the network. This result shows that each node acquires the information contained in the gradient oracles at every other node and suffers asymptotically no regret as the number of gradient oracle calls at each node increse.\\

{\bf Proof }
According to Lemma 1 and Lemma 2, we have
\begin{align}\label{eq:step2}
  & \E \left[ f(\vw_t) - f(\vw^*) \right] \nn\\
    & \hspace{.5cm} \leq \frac{N}{2\mu_t} \, \E \left[ (1-\lambda\mu_t)\norm{\vw_t-\vw^*}^2 - \norm{\vw_{t+1}-\vw^*}^2 \right] \nn\\
    & \hspace{1.5cm} + 3NG^2\mu_t + 6N\sqrt{N}G^2 \sum_{z=1}^{t-1} \mu_{t-z} \sigma^{z}.
\end{align}
From the convexity of the cost functions, we also have
\begin{equation}\label{eq:pot1}
  \E \left[ f_i(\vw_t)-f_i(\vw_{t,j}) \right] \geq \E \, \vg_{t,i,j}^T(\vw_t-\vw_{t,j}),
\end{equation}
$\forall i,j\in\{1,\dots,N\}$, where $\vg_{t,i,j} \triangleq \nabla_{\vw_{t,j}} f_i(\vw_{t,j})$. Here, we can rewrite \eqref{eq:pot1} as follows
\begin{align}\label{eq:pot2}
  \E \left[ f_i(\vw_{t,j})-f_i(\vw_t) \right] & \leq \E \, \vg_{t,i,j}^T(\vw_{t,j}-\vw_t) \nn\\
    & \leq \E \left[ \norm{\vg_{t,i,j}} \norm{\vw_{t,j}-\vw_t} \right] \nn\\
    & \leq G \, \E \norm{\vw_{t,j}-\vw_t},
\end{align}
where the second line follows from the triangle inequality and the last line follows from the boundedness assumption. Summing \eqref{eq:pot2} from $i=1$ to $N$, we obtain
\begin{equation}\label{eq:pot3}
  \E \left[ f(\vw_{t,j})-f(\vw_t) \right] \leq N G \, \E \norm{\vw_{t,j}-\vw_t}.
\end{equation}
Using Lemma 2 in \eqref{eq:pot3}, we get
\begin{equation}\label{eq:pot4}
  \E \left[ f(\vw_{t,j})-f(\vw_t) \right] \leq 2 N\sqrt{N} G^2 \sum_{z=1}^{t-1} \mu_{t-z} \sigma^z.
\end{equation}
We then add \eqref{eq:step2} and \eqref{eq:pot4} to obtain
\begin{align}\label{eq:step3}
  & \E\left[ f(\vw_{t,j}) - f(\vw^*) \right] \nn\\
    & \hspace{.5cm} \leq \frac{N}{2\mu_t} \E\left[ (1-\lambda\mu_t)\norm{\vw_t-\vw^*}^2 - \norm{\vw_{t+1}-\vw^*}^2 \right] \nn\\
    & \hspace{1.5cm} + 3NG^2\mu_t + 8N\sqrt{N}G^2 \sum_{z=1}^{t-1} \mu_{t-z} \sigma^{z}.
\end{align}

Multiplying both sides of \eqref{eq:step2} by $t$ and summing from $t=1$ to $T$ yields \cite{lacoste}
\begin{align}\label{eq:step4}
  & \E \sum_{t=1}^T t \left[ f(\vw_{t,j}) - f(\vw^*) \right] \nn\\
    & \hspace{.1cm} \leq \frac{N(1-\lambda\mu_1)}{2\mu_1} \E\norm{\vw_1-\vw^*}^2 - \frac{TN}{2\mu_T} \E\norm{\vw_{T+1}-\vw^*}^2 \nn\\
    & \hspace{.5cm} + \sum_{t=2}^T \frac{N}{2} \left( \frac{t(1-\lambda\mu_t)}{\mu_t} - \frac{t-1}{\mu_{t-1}} \right) \E\norm{\vw_t-\vw^*}^2 \nn\\
    & \hspace{.5cm} + 3NG^2\sum_{t=1}^T t\mu_t +  8N\sqrt{N}G^2 \sum_{t=1}^T t \sum_{z=1}^{t-1} \mu_{t-z} \sigma^{z}.
\end{align}
Here, we observe that
\begin{align}\label{eq:step4_2}
  \sum_{t=1}^T \sum_{z=1}^{t-1} t \mu_{t-z} \sigma^{z} & \leq \sum_{z=1}^T \sum_{t=1}^T t \mu_{t-z} \sigma^{z} \nn\\
    & \leq \sum_{z=1}^T\sigma^{z}  \sum_{t=1}^T t \mu_t \nn\\
    & \leq \frac{\sigma}{1-\sigma} \sum_{t=1}^T t \mu_t.
\end{align}
Putting \eqref{eq:step4_2} back in \eqref{eq:step4} and inserting $\mu_t=\frac{2}{\lambda(t+1)}$, we obtain
\begin{align}\label{eq:step5}
  & \E \sum_{t=1}^T t \left[ f(\vw_{t,j}) - f(\vw^*) \right]  \nn\\
    & \hspace{.5cm} \leq - \frac{\lambda NT(T+1)}{4} \, \E \norm{\vw_{T+1}-\vw^*}^2 \nn\\
    & \hspace{1.5cm} + \left( 3NG^2 + 8N\sqrt{N}G^2 \frac{\sigma}{1-\sigma} \right) \sum_{t=1}^T \frac{2t}{\lambda(t+1)} \nn\\
    & \hspace{.5cm} \leq - \frac{\lambda NT(T+1)}{4} \, \E \norm{\vw_{T+1}-\vw^*}^2 \nn\\
    & \hspace{1.5cm} + \frac{2NG^2T}{\lambda} \left( 3 + 8\sqrt{N}\frac{\sigma}{1-\sigma} \right),
\end{align}
where the last line follows since $\frac{t}{t+1} \leq 1$. Dividing both sides of \eqref{eq:step5} by $\sum_{t=1}^T t = \frac{T(T+1)}{2}$, we obtain
\begin{align}\label{eq:step6}
  & \E\left[ \frac{2}{T(T+1)} \sum_{t=1}^T t \left( f(\vw_{t,j}) - f(\vw^*) \right) \right] \nn\\
    & \hspace{1.5cm} \leq - \frac{\lambda N}{2} \, \E\norm{\vw_{T+1}-\vw^*}^2 \nn\\
    & \hspace{2.5cm} + \frac{4NG^2}{\lambda(T+1)} \left( 3 + 8\sqrt{N}\frac{\sigma}{1-\sigma} \right).
\end{align}
Since $f_i$'s are convex for all $i\in\{1,\dots,N\}$, $f$ is also convex. Thus, from Jensen's inequality, we can write
\begin{align}\label{eq:step7}
  & \E\left[ f\left( \frac{2}{T(T+1)} \sum_{t=1}^T t \, \vw_{t,j} \right) - f(\vw^*) \right] \nn\\
    & \hspace{2cm} \leq \E\left[ \frac{2}{T(T+1)} \sum_{t=1}^T t \left( f(\vw_{t,j}) - f(\vw^*) \right) \right].
\end{align}
Combining \eqref{eq:step6} and \eqref{eq:step7}, we obtain
\begin{align}\label{eq:step8}
  & \E\left[ f\left( \frac{2}{T(T+1)} \sum_{t=1}^T t \, \vw_{t,j} \right) - f(\vw^*) \right] \nn\\
    & \hspace{1.5cm} \leq - \frac{\lambda N}{2} \, \E\norm{\vw_{T+1}-\vw^*}^2 \nn\\
    & \hspace{2.5cm} + \frac{4NG^2}{\lambda(T+1)} \left( 3 + 8\sqrt{N}\frac{\sigma}{1-\sigma} \right).
\end{align}
This concludes the proof of Theorem 1. \hfill $\square$\\

Hence, using the weighted average
\begin{equation}\label{eq:tvw}
  \vbw_{T,j} \triangleq \frac{2}{T(T+1)} \sum_{t=1}^T t \, \vw_{t,j},
\end{equation}
instead of the original SGD iterates, we can achieve a convergence rate of $O\left(\frac{N\sqrt{N}}{T}\right)$. The denominator $T$ of this regret bound follows since we use a time variable weighting of the SGD iterates. The linear dependency to the network size follows since we add $N$ different cost functions, i.e., one corresponding to each node. Finally, the sublinear dependency to the network size results from the diffusion of the parameter vector over the distributed network.

In the following theorem, we then consider the performance of the average SGD iterate instead of the time-variable weighted iterate in \eqref{eq:tvw}. In particular, we show that the average SGD iterate achieves a MSD of $O\left(\frac{\sqrt{N}}{T}\right)$. This MSD follows due to the number of gradient oracle calls and diffusion regret over the distributed network.\\

{\bf Theorem 2}
{\em Under the assumptions in Lemma 1 and Lemma 2, Algorithm \ref{alg} with learning rate $\mu_t=\frac{2}{\lambda(t+1)}$ and weighted parameters $\vbw_{t,i}$, when applied to any independent and identically distributed regressor vectors and data, i.e., $(\vu_{t,i},d_{t,i})$ for all $t\geq1$ and $i=1,\dots,N$, yields the following MSD guarantee
\begin{equation}\label{eq:theorem2}
  \E \norm{\vw_{T+1}-\vw^*}^2 \leq \frac{24G^2}{\lambda^2(T+1)} \left( 1 + \frac{2\sigma\sqrt{N}}{1-\sigma} \right).
\end{equation}
for all $T \geq 1$, where $\vw_t = \frac{1}{N}\sum_{i=1}^N \vw_{t,i}$ and $\sigma$ is the second largest singular value of matrix $\vH$.\\
}

{\bf Proof }
Following similar lines to the proof of Theorem 1, in particular following the steps between \eqref{eq:step3}-\eqref{eq:step8}, we get
\begin{align}\label{eq:step9}
  & \E\left[ f\left( \frac{2}{T(T+1)} \sum_{t=1}^T t \, \vw_t \right) - f(\vw^*) \right] \nn\\
    & \hspace{1.5cm} \leq - \frac{\lambda N}{2} \, \E\norm{\vw_{T+1}-\vw^*}^2 \nn\\
    & \hspace{2.5cm} + \frac{12NG^2}{\lambda(T+1)} \left( 1 + 2\sqrt{N}\frac{\sigma}{1-\sigma} \right),
\end{align}
which yields
\begin{equation}\label{eq:step10}
  \E\norm{\vw_{T+1}-\vw^*}^2 \leq \frac{24G^2}{\lambda^2(T+1)} \left( 1 + 2\sqrt{N}\frac{\sigma}{1-\sigma} \right),
\end{equation}
which is the desired result. This concludes the proof of Theorem 2. \hfill $\square$\\

%
%
%
%
%

\section{Extensions to the Consensus Strategy}\label{sec:ext}
Algorithm \ref{alg} can be generalized into the consensus strategy in a straightforward manner, while the performance guarantee in Theorem 1 still holds up to constant terms, i.e., we still have a convergence rate of $O\left(\frac{N\sqrt{N}}{T}\right)$. For the consensus strategy, the lines 5--7 of Algorithm \ref{alg} is replaced by the following update
\begin{equation}
  \vw_{t+1,i} = \Pi_\W\left( \sum_{j=1}^N \vH_{ji} \vw_{t,j} - \mu_t \vg_{t,i} \right).
\end{equation}
Hence, we have the following recursion on the parameter vectors
\begin{equation}
  \vW_t = \vW_1 \vH^{t-1} - \sum_{z=1}^{t-1} \left( \vX_{t-z} - \mu_{t-z} \vG_{t-z} \right) \vH^{z-1},\nn
\end{equation}
instead of the one in \eqref{eq:recursion}. According to this modification, Lemma 2 can be updated as follows
\begin{equation}
  \E\norm{\vw_t-\vw_{t,i}} \leq 2 G \sqrt{N} \sum_{z=1}^{t-1} \mu_{t-z} \sigma^{z-1}.
\end{equation}
This loosens the upper bounds in \eqref{eq:bound4} and \eqref{eq:bound5} by a factor of $1/\sigma$. Therefore, diffusion strategies achieves a better convergence performance compared to the consensus strategy.

\section{Simulations}\label{sec:sim}
\begin{figure}
  \centering
  \includegraphics[width=.48\textwidth]{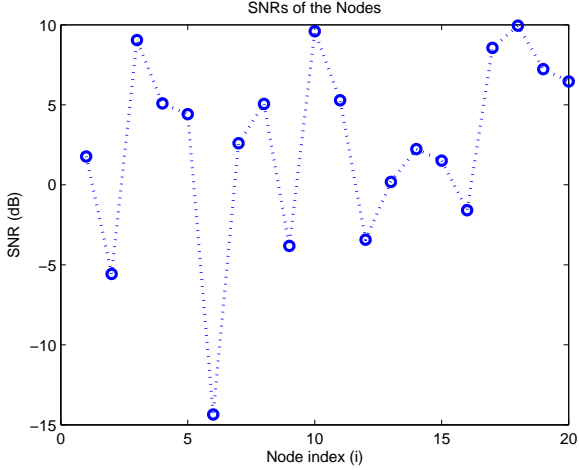}\\
  \caption{SNRs of the nodes in the distributed network.}\label{fig:snr}
\end{figure}

\begin{figure*}[!t]
    \centering
    \subfloat[Global NCE for the star network]{\begin{overpic}[width=0.48\textwidth]{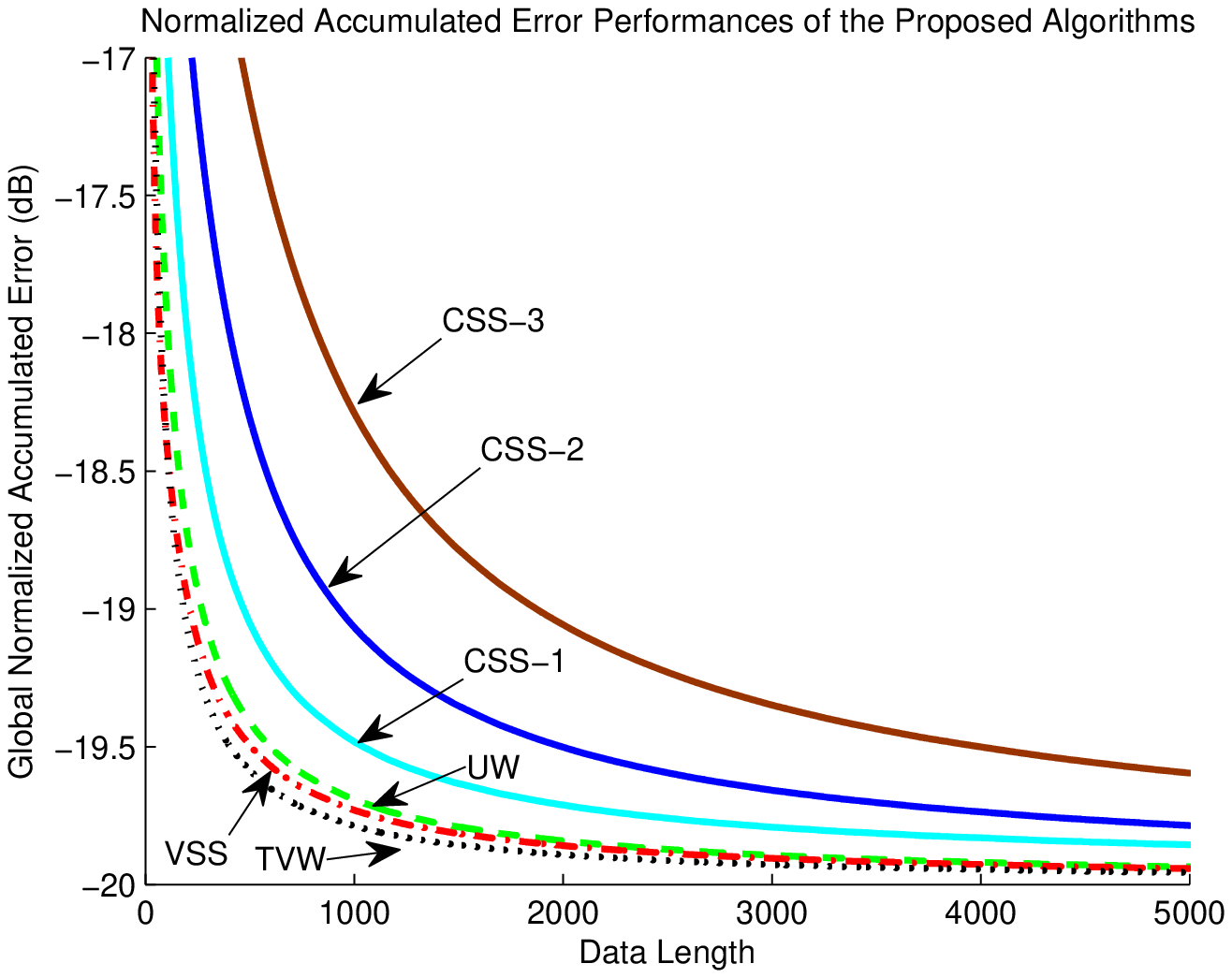}
     \put(55,40){\includegraphics[scale=0.2]{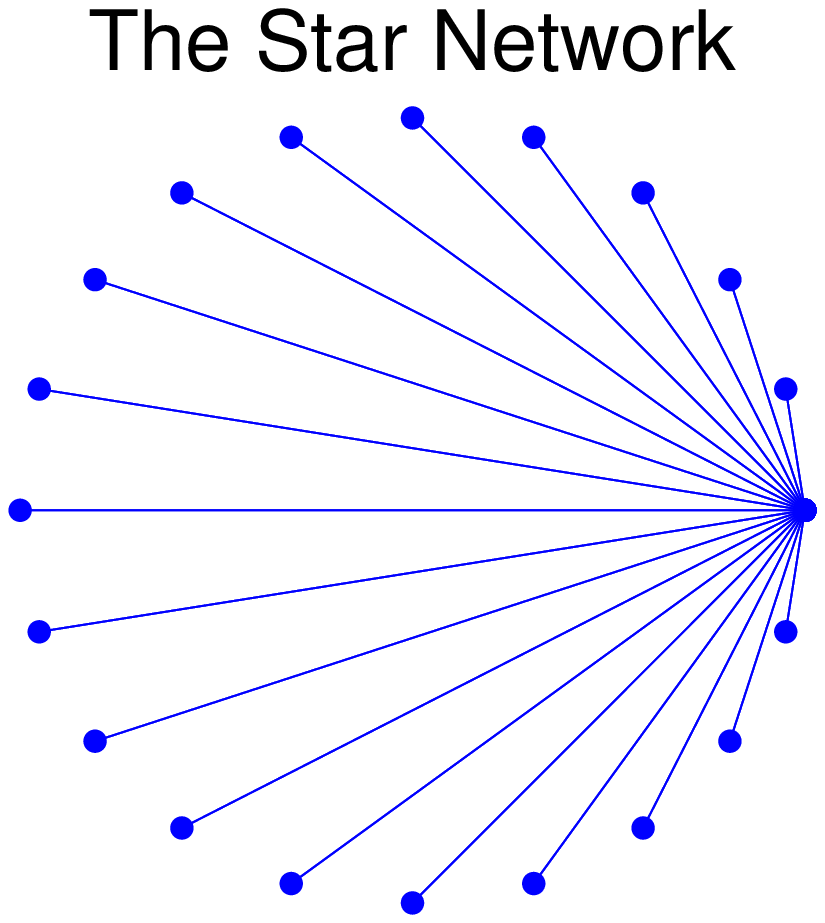}}
  \end{overpic}}
    \hfil
    \subfloat[Global MSD for the star network]{\includegraphics[width=0.48\textwidth]{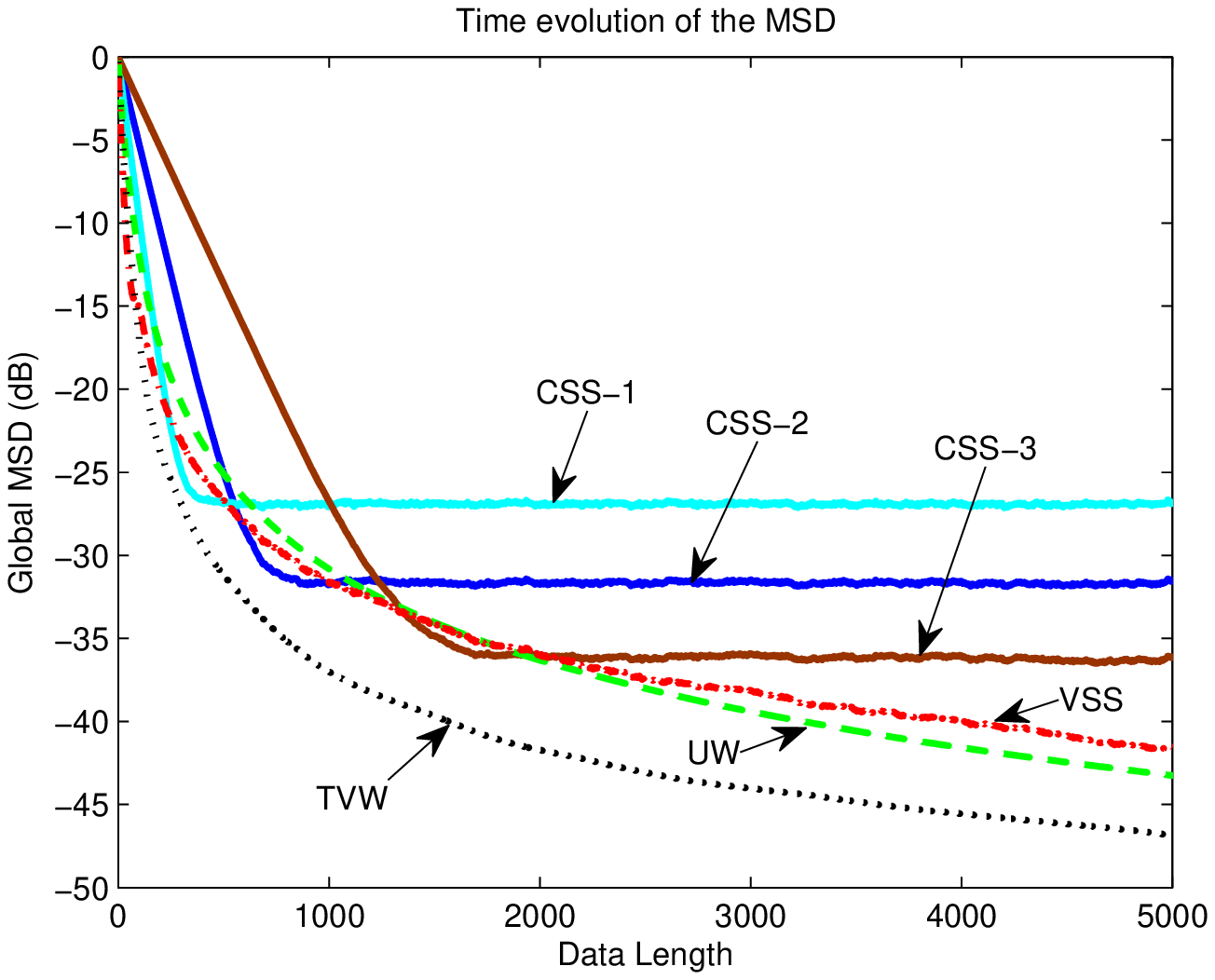}}
    \hfil
    \subfloat[Global NCE for the circle network]{\begin{overpic}[width=0.48\textwidth]{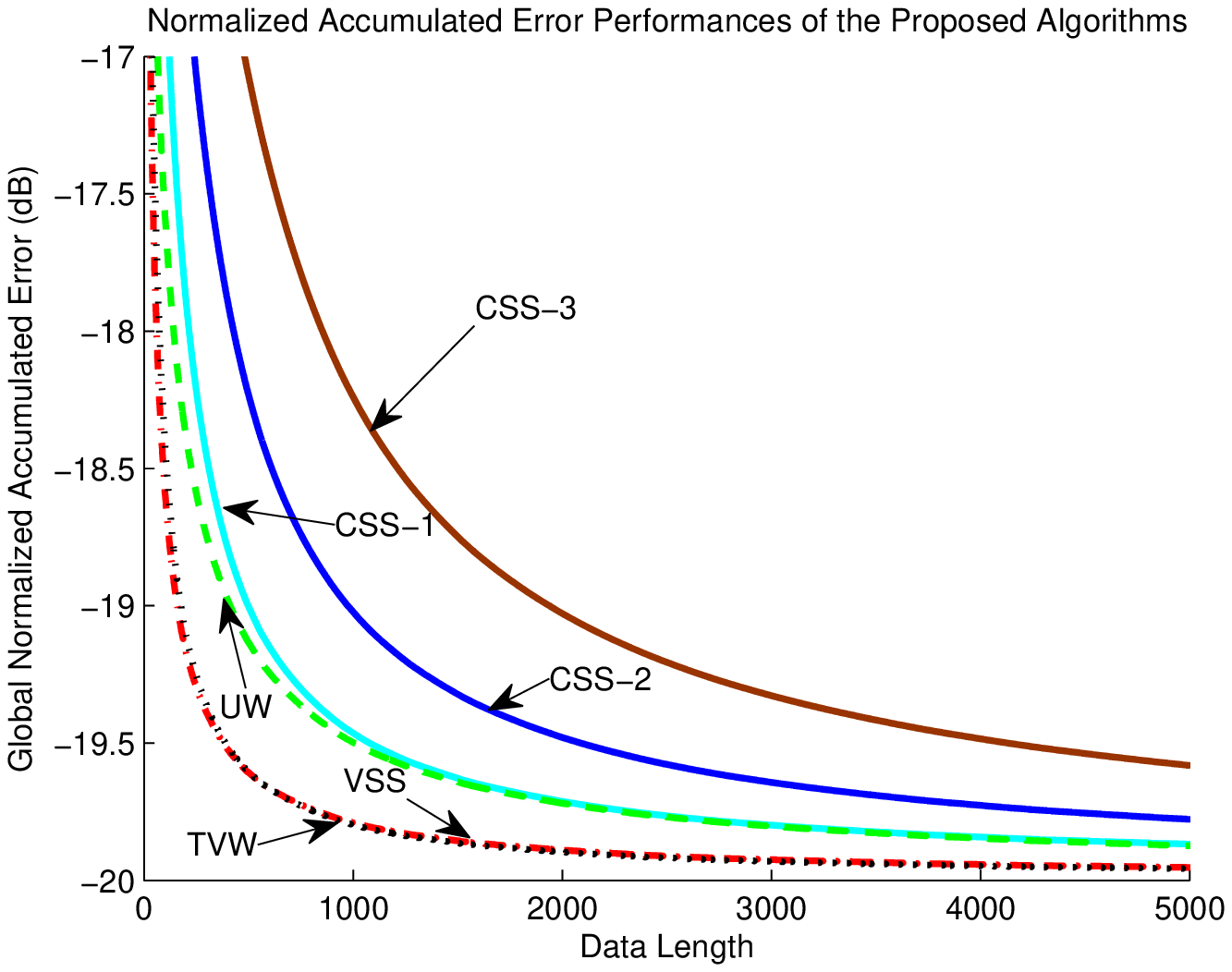}
     \put(55,40){\includegraphics[scale=0.2]{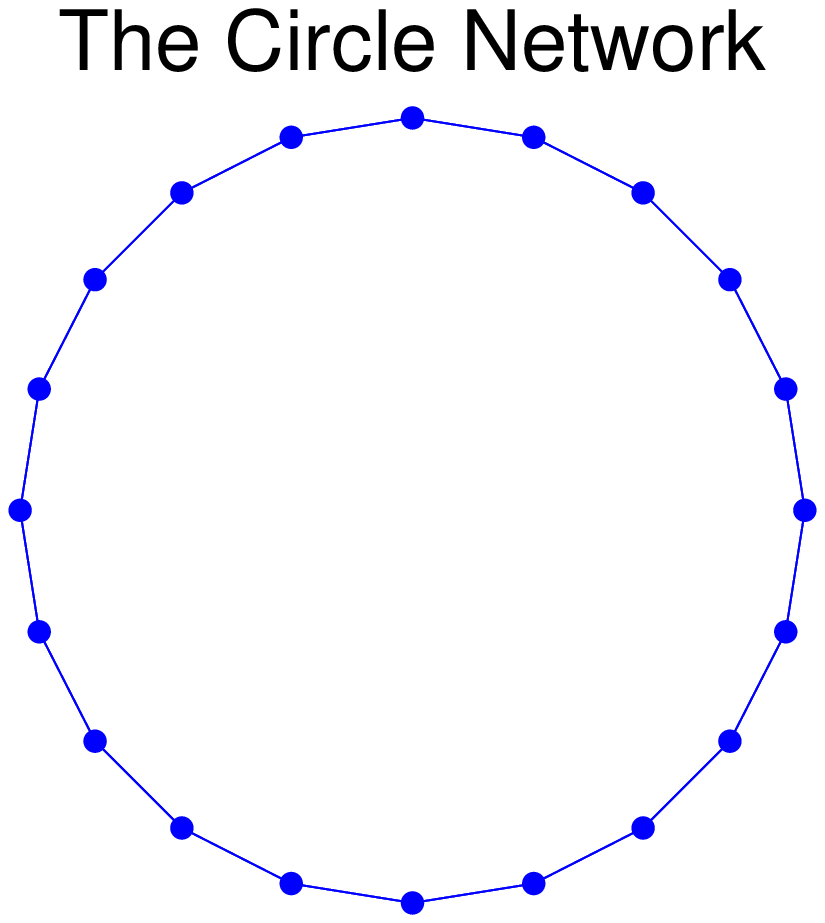}}
  \end{overpic}}
    \hfil
    \subfloat[Global MSD for the circle network]{\includegraphics[width=0.48\textwidth]{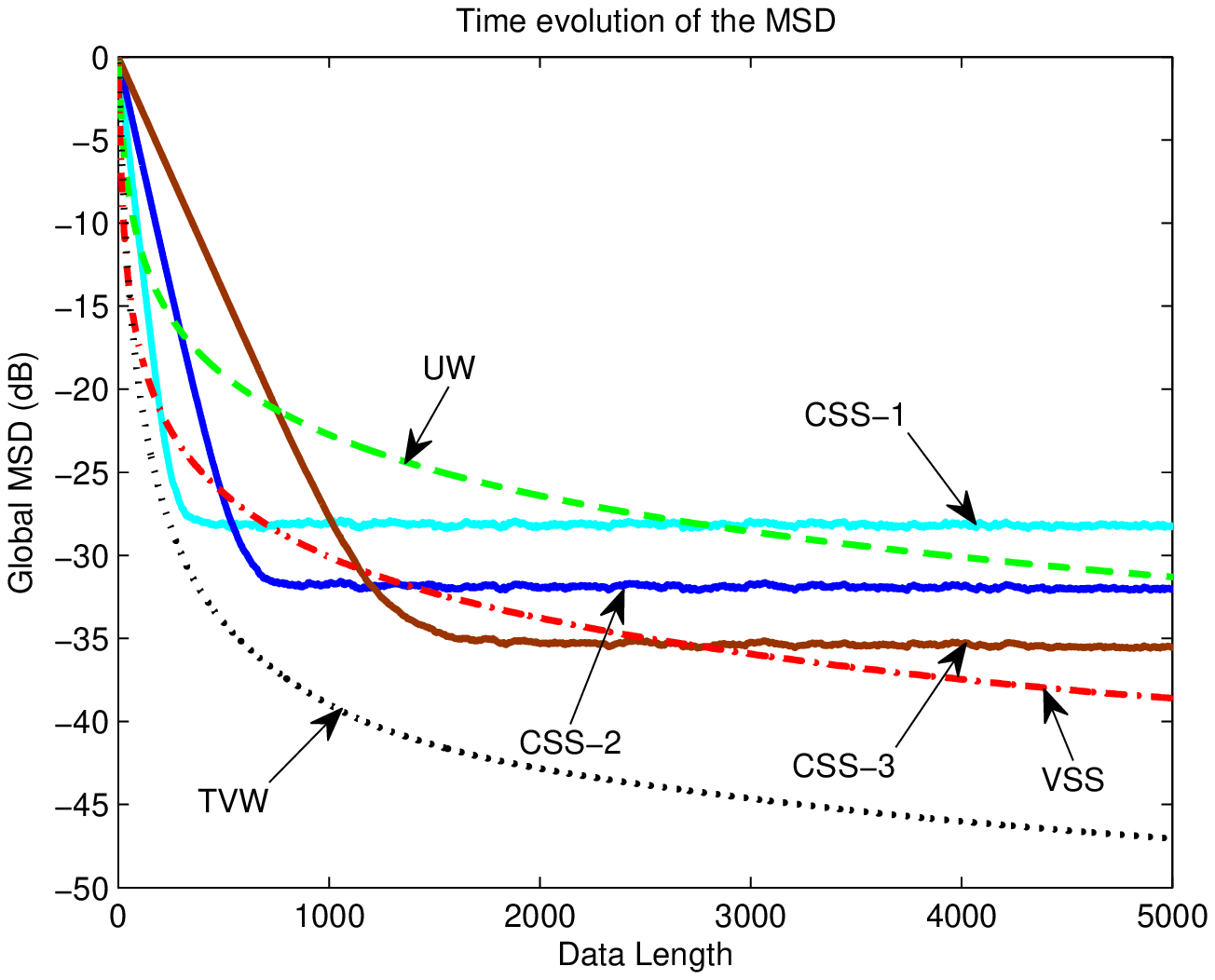}}
    \hfil
    \subfloat[Global NCE for the random network]{\begin{overpic}[width=0.48\textwidth]{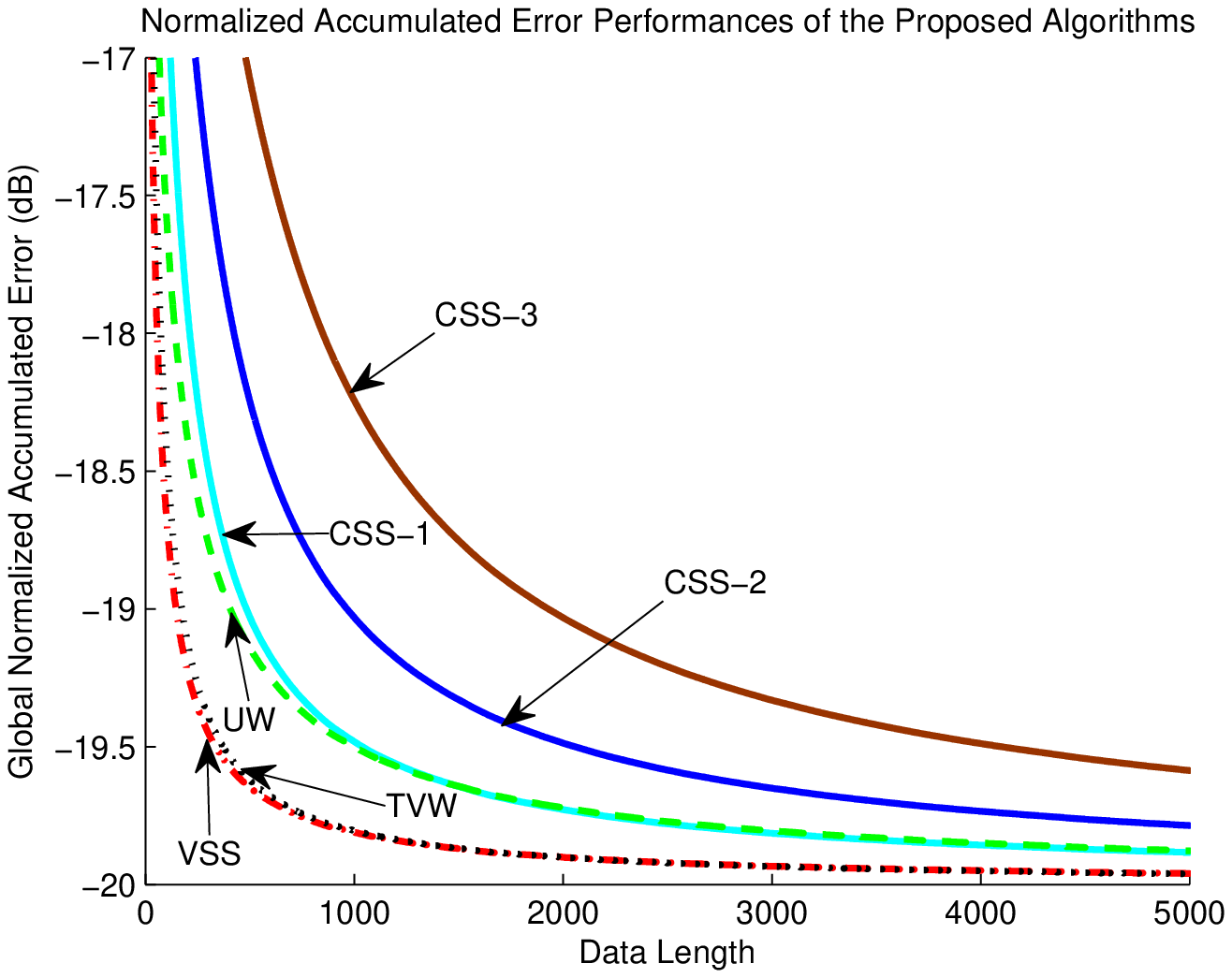}
     \put(55,40){\includegraphics[scale=0.22]{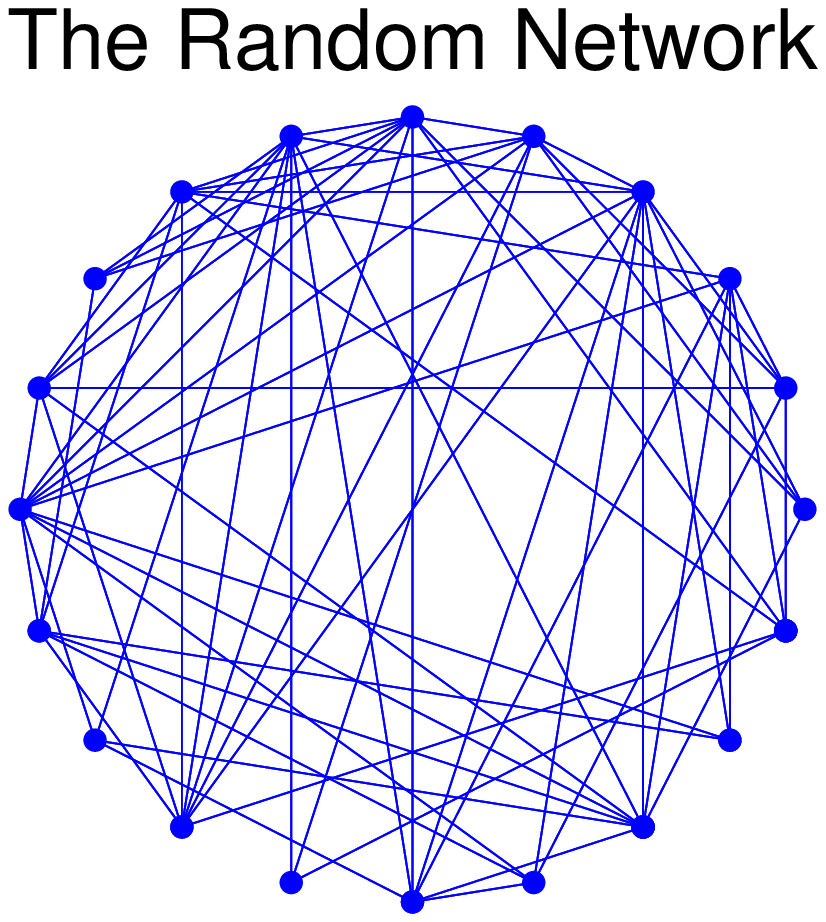}}
  \end{overpic}}
    \hfil
    \subfloat[Global MSD for the random network]{\includegraphics[width=0.48\textwidth]{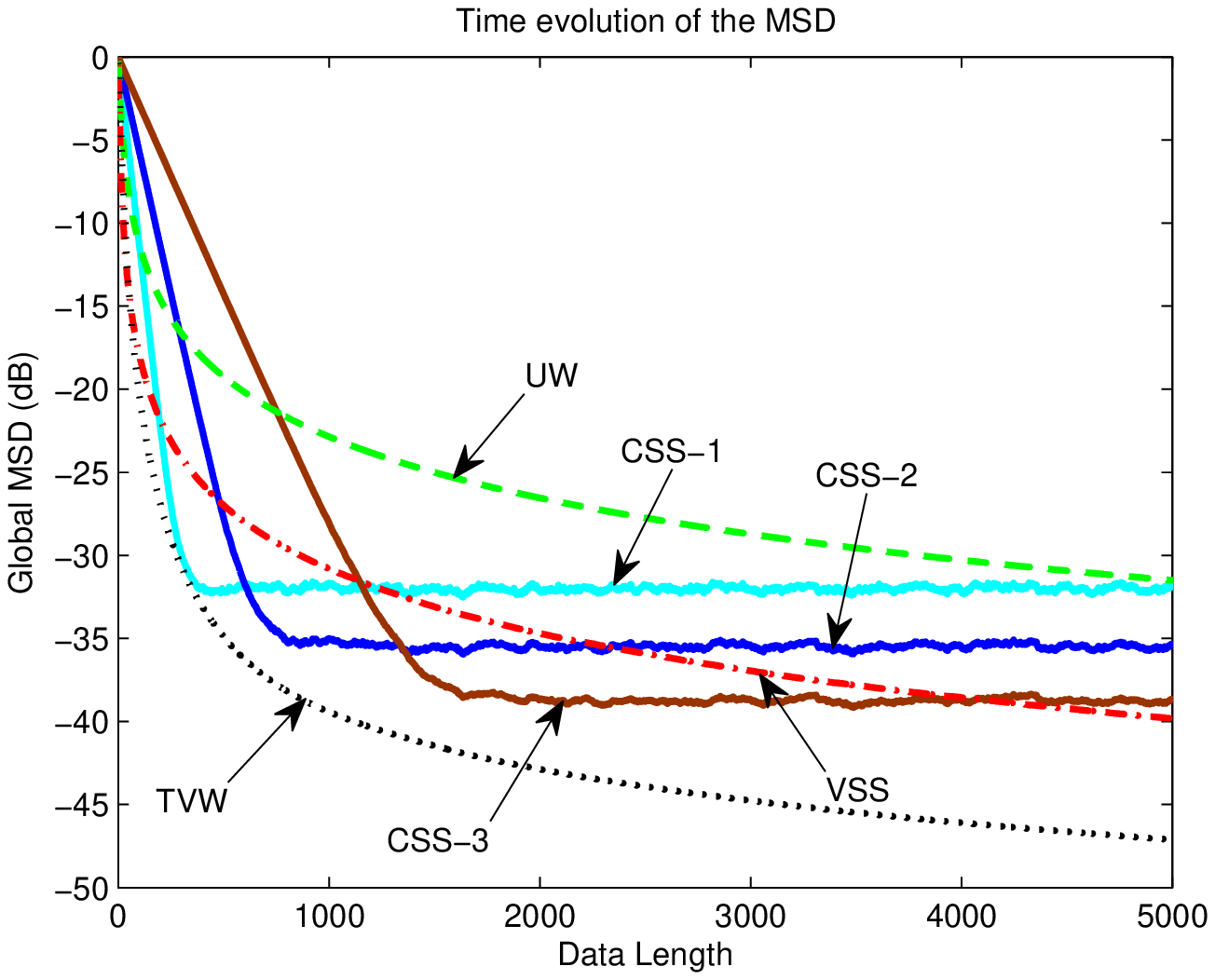}}
    \caption{NCE (left column) and MSD (right column) performances of the proposed algorithms under the star (first row), the circle (second row), and a random (third row) network topologies, under the squared error loss function averaged over $200$ trials.}
    \label{fig}
\end{figure*}

In this section, we first examine the performance of the proposed algorithms for various distributed network topologies, namely the star, the circle, and a random network topologies (which are shown in Fig. \ref{fig}). In all cases, we have a network of $N = 20$ nodes where at each node $i=1,\dots,N$ at time $t$, we observe the data $d_{t,i} = \vec{w}_0^T\vu_{t,i} + v_{t,i}$, where the regression vector $\vu_{t,i}$ and the observation noise $v_{t,i}$ are generated from i.i.d. zero mean Gaussian processes for all $t\geq1$. The variance of the observation noise is $\sigma_{v,i}^2 = 0.1$ for all $i=1,\dots,N$, whereas the auto-covariance matrix of the regression vector $\vu_{t,i} \in \Real^5$ is randomly chosen for each node $i=1,\dots,N$ such that the signal-to-noise ratio (SNR) over the network varies between $-15$dB to $10$dB (see Fig. \ref{fig:snr}). The parameter of interest $\vec{w}_0 \in \Real^5$ is randomly chosen from a zero mean Gaussian process and normalized to have a unit norm, i.e., $\norm{\vec{w}_0} = 1$. We use the well-known Metropolis combination rule \cite{sayed_mag1} to set the combination weights as follows
\begin{equation}
  \vH_{ij} =
  \begin{cases}
    \frac{1}{\max\left\{n_i,n_j\right\}} & \text{, if $j \in {\cal N}_i \setminus i$}\\
    0 & \text{, if $j \notin {\cal N}_i$} \\
    1-\sum_{j \in {\cal N}_i \setminus i} \vH_{ij} & \text{, if $i=j$}
  \end{cases}
\end{equation}
where $n_i$ is the number of neighboring nodes for node $i$.

For this set of experiments, we consider the squared error loss, i.e., $\ell(\vw_{t,i};\vu_{t,i},d_{t,i}) = (d_{t,i} - \vw_{t,i}^T\vu_{t,i})^2$ as our loss function. In the figures, CSS represents the distributed constant step-size SGD algorithm of \cite{sayed2}, VSS represents the distributed variable step-size SGD algorithm of \cite{regret}, UW represents the distributed version of the uniform weighted SGD algorithm of \cite{rakhlin}, and TVW represents the distributed time variant weighted SGD algorithm proposed in this paper. The step-sizes of the CSS-1, CSS-2, and CSS-3 algorithms are set to $0.05$, $0.1$, and $0.2$, respectively, at each node and the learning rates of the VSS and UW algorithms are set to $1/(\lambda t)$ as noted in \cite{regret,rakhlin}, whereas the learning rate of the TVW algorithm is set to $2/(\lambda(t+1))$ as noted in Theorem 1, where $\lambda=0.01$. These learning rates chosen to guarantee a fair performance comparison between these algorithms according to the corresponding algorithm descriptions stated in this paper and in \cite{regret,rakhlin}.

In the left column of Fig. \ref{fig}, we compare the normalized time accumulated error performances of these algorithms under different network topologies in terms of the global normalized cumulative error (NCE) measure, i.e.,
\begin{equation}
  \mathrm{NCE}(t) = \frac{1}{Nt} \sum_{i=1}^N\sum_{\tau=1}^{t} (d_{\tau,i} - \vw_{\tau,i}^T\vu_{\tau,i})^2. \nn
\end{equation}
Additionally, in the right column of Fig. \ref{fig}, we compare the performance of the algorithms in terms of the global MSD measure, i.e.,
\begin{equation}
  \mathrm{MSD}(t) = \frac{1}{N} \sum_{i=1}^N \norm{\vec{w}_0-\vw_{t,i}}^2. \nn
\end{equation}
In the figures, we have plotted the NCE and MSE performances of the proposed algorithms over $200$ independent trials to avoid any bias.

As can be seen in the Fig. \ref{fig}, the proposed TVW algorithm significantly outperforms its competitors and achieves a smaller error performance. This superior performance of our algorithm is obtained thanks to the time-dependent weighting of the regression parameters, which is used to obtain a faster convergence rate with respect to the rest of the algorithms. Hence, by using a certain time varying weighting of the SGD iterates, we obtain a significantly improved convergence performance compared to the state-of-the-art approaches in the literature. Furthermore, the performance of our algorithm is robust against the network topology, whereas the competitor algorithms may not provide satisfactory performances under different network topologies.

We next consider the classification tasks over the benchmark data sets: covertype\footnote{https://www.csie.ntu.edu.tw/$\sim$cjlin/libsvmtools/datasets/} and quantum\footnote{http://osmot.cs.cornell.edu/kddcup/}. For this set of experiments, we consider the hinge loss, i.e., $\ell(\vw_{t,i};\vu_{t,i},d_{t,i}) = \max\{0,1-d_{t,i}\vw_{t,i}^T\vu_{t,i}\}^2$ as our loss function. The regularization constant is set to $\lambda=1/T$, where the stepsizes of the TVW, UW, and VSS algorithms are set as in the previous experiment. The step-sizes of the CSS-1, CSS-2, and CSS-3 algorithms are set to $0.02$, $0.05$, and $0.1$ for the covertype data set, whereas the step-sizes of the CSS-1, CSS-2, and CSS-3 algorithms are set to $0.01$, $0.02$, and $0.05$ for the quantum data set. These learning rates are chosen to illustrate the tradeoff between the convergence speed and the steady state performance of the constant stepsize SGD methods. The network sizes are set to $N=20$ and $N=50$ for the covertype and quantum data sets, respectively.

\begin{figure}
  \centering
  \includegraphics[width=.48\textwidth]{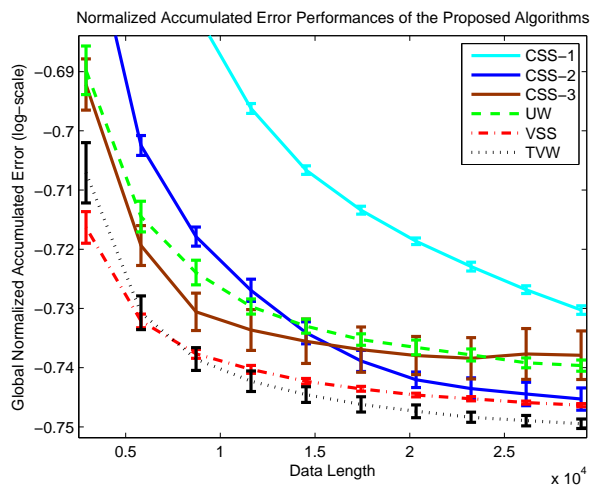}\\
  \caption{Normalized accumulated errors of the proposed algorithms versus training data length for cover type data averaged over $250$ trials for a network size of $20$.}\label{fig:covtype}
\end{figure}

\begin{figure}
  \centering
  \includegraphics[width=.48\textwidth]{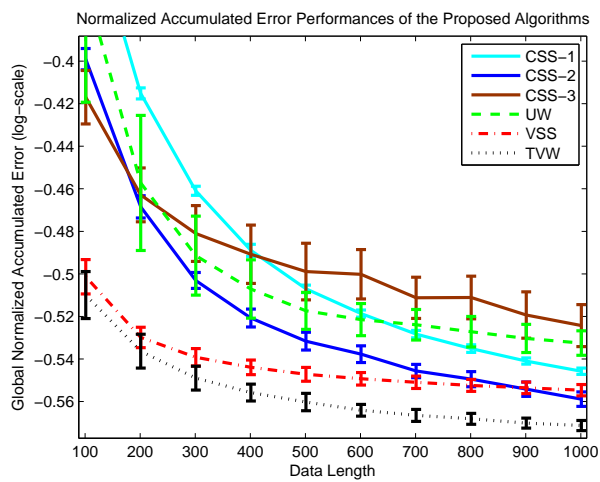}\\
  \caption{Normalized accumulated errors of the proposed algorithms versus training data length for quantum data averaged over $100$ trials for a network size of $50$.}\label{fig:quantum}
\end{figure}

In Fig. \ref{fig:covtype} and Fig. \ref{fig:quantum}, we illustrate the performances of the proposed algorithms for various training data lengths. In particular, we train the parameter vectors at each node using a certain length of training data and test the performance of the final parameter vector over the entire data set. We provide averaged results over $250$ and $100$ independent trials for covertype and quantum data sets, respectively, and present the mean and variance of the normalized accumulated hinge errors. These figures illustrate that the proposed TVW algorithm significantly outperforms its competitors. Although the performances of the UW and VSS algorithms are comparably robust over different iterations, the TVW algorithm provides a smaller accumulated loss. On the other hand, the variance of the constant stepsize methods highly deteriorate as the stepsize increases. Although decreasing the stepsize yields more robust performance for these constant stepsize algorithms, the TVW algorithm provides a significantly smaller steady-state cumulative error with respect to these methods.

\section{Conclusion}\label{sec:conc}
We study distributed strongly convex optimization over distributed networks, where the aim is to minimize a sum of unknown convex objective functions. We introduce an algorithm that uses a limited number of gradient oracle calls to these objective functions and achieves the optimal convergence rate of $O\left(\frac{N\sqrt{N}}{T}\right)$ after $T$ gradient updates at each node. This performance is obtained by using a certain time-dependent weighting of the SGD iterates at each node. The computational complexity and the communication load of the proposed approach is the same with the state-of-the-art methods in the literature up to constant terms. We also prove that the average SGD iterate achieves a mean square deviation (MSD) of $O\left(\frac{\sqrt{N}}{T}\right)$ after $T$ gradient oracle calls. We illustrate the superior convergence rate of our algorithm with respect to the state-of-the-art methods in the literature.


\end{document}